\documentclass[%
 aip,
 amsmath,amssymb,
preprint,%
]{revtex4-2}

\usepackage{booktabs}

\usepackage{graphicx}
\usepackage{dcolumn}
\usepackage{bm}

\usepackage[utf8]{inputenc}
\usepackage[T1]{fontenc}
\usepackage{etoolbox}

\graphicspath{{./eps/}{./epstopdf_converted/}}

\newcommand*\D {\mathcal{D}}
\newcommand*\F {\mathcal{F}}

\newcommand* \disp {\phi}

\newcommand*\z {\varepsilon}

\newcommand*\kk {\kappa}

\newcommand* \RevisionDelete[1] {}

\usepackage{xcolor}

\makeatletter
\def\@email#1#2{%
 \endgroup
 \patchcmd{\titleblock@produce}
  {\frontmatter@RRAPformat}
  {\frontmatter@RRAPformat{\produce@RRAP{*#1\href{mailto:#2}{#2}}}\frontmatter@RRAPformat}
  {}{}
}%
\makeatother

\begin{document}



	\title{ 
    Koopman Theory-Inspired Method for Learning Time Advancement Operators in Unstable Flame Front Evolution }
	\author{Rixin Yu}
	\email{rixin.yu@energy.lth.se}
	\affiliation{ 
	Department of Energy Sciences, Lund University, 22100 Lund, Sweden
	}%
	\author{Marco Herbert}
	\affiliation{Department of Aerospace Engineering, 
University of the Bundeswehr Munich, Germany}%
	\author{Markus Klein}
	\affiliation{Department of Aerospace Engineering, 
University of the Bundeswehr Munich, Germany}%
	\author{Erdzan Hodzic}
\affiliation{%
Department of Manufacturing Processes, RISE Research Institutes of Sweden, 553 22 Jonkoping, Sweden 
}%

	\date{\today}
	
	\begin{abstract}
		Predicting the evolution of complex systems governed by partial differential equations (PDEs) remains challenging, especially for nonlinear, chaotic behaviors. This study introduces Koopman-inspired Fourier Neural Operators (kFNO) and Convolutional Neural Networks (kCNN) to learn solution advancement operators for flame front instabilities. By transforming data into a high-dimensional latent space, these models achieve more accurate multi-step predictions compared to traditional methods. Benchmarking across one- and two-dimensional flame front scenarios demonstrates the proposed approaches' superior performance in short-term accuracy and long-term statistical reproduction, offering a promising framework for modeling complex dynamical systems.
	\end{abstract}
	
	\maketitle



\section{Introduction \label{sec:intro}} 

Partial differential equations (PDEs) are fundamental mathematical frameworks used to describe complex physical phenomena across diverse scientific and engineering domains. From fluid dynamics and climate modeling to quantum mechanics and biological systems, PDEs encapsulate intricate interactions and dynamical behaviors derived from underlying physical principles. However, solving PDEs, particularly nonlinear equations with complex boundary conditions, poses significant computational challenges, historically limiting our ability to simulate and predict such systems accurately.

The computational landscape for solving PDEs has been transformed by the integration of machine learning (ML) and artificial intelligence (AI) techniques. Recent advancements have introduced a proliferation of operator learning methods, each contributing unique insights and capabilities for tackling complex mathematical problems. Early efforts in this domain leveraged convolutional neural networks (CNNs)\cite{CNN1,CNN2,CNN3,CNN4,CNN5,UNet,ConvPDE}, drawing inspiration from computer vision to parameterize PDE operators in finite-dimensional spaces. These CNN-based approaches mapped discrete functions to image-like representations, offering innovative computational frameworks for PDE solutions.

Building on CNN-based methods, neural operator techniques \cite{GraphKenerlNetwork, kovachki2021neural}  emerged to learn operators of infinite dimensionality, achieving remarkable success across diverse benchmarks. Notable examples include Deep Operator Networks (DeepOnet)\cite{DeepONet} and Fourier Neural Operators (FNO)\cite{FNO2020}. These methods introduced sophisticated strategies for mapping between function spaces: FNO leverages Fourier transforms for efficient learning of continuous operators, while DeepOnet employs a branch-trunk network architecture to model complex relationships. Recent advancements have further extended these approaches, incorporating wavelet-based methods \cite{gupta2021multiwavelet, tripura2023wavelet}, adaptations to complex domains \cite{chen2023laplace}, and techniques for parametric-dependent operator learning \cite{YuH2024}.

Among the various operators encountered in PDE problems, the solution time-advancement operator is of particular interest. This operator maps the solution at a given time to the solution at a future time, over a short interval, under prescribed boundary conditions. By recurrently applying this operator — using its output as the input for the next step — it becomes possible to generate predictions over long time horizons. Applications such as weather forecasting, transient fluid flow simulations, and flame front evolution rely heavily on such operator learning approaches to capture dynamic evolution across spatial and temporal scales.

Learning these operators presents unique challenges. While accurate short-term predictions are crucial, ensuring stability and fidelity in long-term predictions becomes particularly difficult for PDEs exhibiting chaotic solutions, where errors grow exponentially due to sensitivity to initial conditions. Previous studies have tackled this issue using recurrent training strategies\cite{FNO2020,Yu2023}, minimizing multi-step errors over consecutive predictions. By contrast, non-recurrent training, which optimizes for single-step prediction accuracy, often leads to divergence in long-term solutions, especially for chaotic systems\cite{Yu2023}.

Inspired by Koopman operator theory\cite{koopman1931hamiltonian,mezic2013analysis,brunton2022data}, we propose a novel approach to directly learn an extended solution advancement operator. The Koopman operator is a linear operator defined in an infinite-dimensional functional space of system observables, offering a powerful analytical framework for nonlinear dynamics. In this work, we adapt this concept to develop Koopman-inspired Fourier Neural Operators (kFNO) and Koopman-inspired Convolutional Neural Networks (kCNN). These architectures build on our previous findings\cite{Yu2023,YuH2024, YHN2024}, which identified that two baseline approaches of FNO and CNN are suitable for learning time-advancement operators on chaotic PDEs solutions.

Our approach transforms data into a high-dimensional latent space, where it is advanced iteratively to capture multi-step dynamics. This latent space transformation is a direct extension of existing operator learning methods, which lift low-dimensional input data into higher-dimensional spaces before projecting it back. To demonstrate the effectiveness of the proposed approach, we apply it to a canonical problem in combustion science: the evolution of unstable flame fronts driven by two intrinsic instabilities — the Darrieus-Landau (DL) and Diffusive-Thermal (DT) mechanisms. The former arises from hydrodynamic effects across density gradients\cite{DARRIEUS1938UNPB,landau1988theory}, while the latter stems from disparities in heat and reactant diffusion\cite{zeldovich1944selected,sivashinsky1977diffusional}. These instabilities are modeled using two nonlinear PDEs: the Michelson-Sivashinsky(MS) \cite{michelson1977nonlinear} and Kuramoto-Sivashinsky(KS) \cite{kuramoto1978diffusion} equations, respectively.

We benchmark our kFNO and kCNN models against baseline methods in both one- and two-dimensional configurations. Metrics for comparison include short-term prediction accuracy and long-term statistics in reproducing nonlinear front evolution, especially under chaotic scenarios.

The paper is organized as follows: Section \ref{sec:OLmethods} details the problem setup for learning solution advancement operators and section \ref{sec:networks} describes the proposed Koopman-inspired operator learning methods. Section \ref{sec:results} shows the results on benchmarking the proposed methods against baseline approaches on 1D and 2D versions of the MS and KS equations. Finally, Section \ref{sec:conclusion} summarizes the findings.

\section{Problem Setup for Operator Learning  \label{sec:OLmethods} }

In this section, we outline the problem setup for learning a PDE operator, alongside a description of recurrent training methods.

Consider a system governed by a partial differential equation (PDE), typically involving multiple functions and mappings between them. A general operator mapping could be denoted as
\begin{equation}
    G : \mathcal{V} \to \mathcal{V}'; \quad v(x) \mapsto v'(x'),
\end{equation}
where the input function $v(x)$, with $x \in \mathcal{D}$, belongs to the functional space $\mathcal{V}(\mathcal{D}, \mathbb{R}^{d_v})$, with domain $\mathcal{D} \subset \mathbb{R}^d$ and codomain $\mathbb{R}^{d_v}$. The output function $v'(x')$, where $x' \in \mathcal{D}'$, resides in another functional space $\mathcal{V}'(\mathcal{D}', \mathbb{R}^{d_v'})$, with domain $\mathcal{D}' \subset \mathbb{R}^{d'}$ and codomain $\mathbb{R}^{d_v'}$.

In this work we are primarily interested in the time advancement operator for the solution of the PDE, defined as
\begin{equation}
    G: \disp(x, t) \mapsto \disp(x, t_1),
    \label{eq:G_operator}
\end{equation}
where $\disp(x, t)$ is the solution to the PDE, and $t_j = t + j\Delta_t$ for $j = 1, 2, \dots$ represents a future time sequence separated by a small time increment $\Delta_t$. For simplicity, we assume that the input and output functions share the same domain and codomain, i.e., $\mathcal{D}' = \mathcal{D}$, $\mathcal{V}' = \mathcal{V}$, $d' = d$, and $d_v' = d_v$, with periodic boundary conditions on $\mathcal{D}$.

In some cases, it may be preferable to work with a related operator that outputs solutions over multiple consecutive time steps:
\begin{equation}
    \bar{G}: \mathcal{V} \to \mathcal{V}^n; \quad \disp(x, t) \mapsto (\disp(x, t_1), \disp(x, t_2), \dots, \disp(x, t_n)),
    \label{eq:Gn_operator}
\end{equation}
where $\mathcal{V}^n$ denotes the Cartesian product of $n$ copies of $\mathcal{V}$. These two operators are related as
\begin{equation}
    \bar{G} = (G, G^2, \dots, G^n),
    \label{eq:Gn_G}
\end{equation}
where the superscript $n$ over an operator denotes its $n$-th power, i.e., iterative application of $G$ for $n$ times: $G^n = G \circ \dots \circ G$.

To approximate the operators $G$ and $\bar{G}$ using neural network methods, we define the networks as follows:
\begin{equation}
    \mathcal{G}: \mathcal{V} \times \Theta \to \mathcal{V}, \quad \text{or equivalently}, \quad \mathcal{G}_{\theta}: \mathcal{V} \to \mathcal{V}, \quad \theta \in \Theta,
\end{equation}
and
\begin{equation}
    \bar{\mathcal{G}}: \mathcal{V} \times \bar{\Theta} \to \mathcal{V}^n, \quad \text{or equivalently}, \quad \bar{\mathcal{G}}_{\bar{\theta}}: \mathcal{V} \to \mathcal{V}^n, \quad \bar{\theta} \in \bar{\Theta}.
\end{equation}
Here, $\Theta$ and $\bar{\Theta}$ denote the spaces of trainable parameters in the respective networks. The training objective is to find the optimal parameters $\theta^* \in \Theta$ (or $\bar{\theta}^* \in \bar{\Theta}$) such that $\mathcal{G}_{\theta^*}$ approximates $G$ (or $\bar{\mathcal{G}}_{\bar{\theta}^*}$ approximates $\bar{G}$).

Our prior studies \cite{Yu2023, YuH2024,YHN2024} focused on learning the single-step time advancement operator $G$. Starting with an initial solution function $\disp(x, t)$, recurrent application of the learned operator $\mathcal{G}_{\theta}$ allows for rolling out predictions over an arbitrary time horizon by iteratively updating the input function with the output from the previous prediction. An ideally learned operator should yield accurate short-term predictions. While long-term predictions may deviate if the underlying PDE has chaotic solutions, it is still desirable for long-term predictions to capture the correct statistical behavior.

To achieve this goal, the training approach in \cite{Yu2023, YuH2024,YHN2024} employed a one-to-many setup, where the network $\mathcal{G}_{\theta}$ was trained to make multiple successive predictions from a single input. This setup ensures the learned solution advancement operator is numerically stable, a crucial consideration highlighted in prior works \cite{Yu2023, YuH2024,YHN2024}. Specifically, given a training set 
\[
    \left\{v_j, (G^1 v_j, \dots, G^n v_j)\right\}_{j=1}^Z,
\]
composed of $Z$ input-output pairs arranged in a 1-to-$n$ manner, training the network $\mathcal{G}_{\theta}$ to approximate $G$ becomes a minimization task:
\begin{equation}
    \min_{\theta \in \Theta} \mathbb{E}_{v \sim \chi} \left[C\left((\mathcal{G}_{\theta}^1 v, \dots, \mathcal{G}_{\theta}^n v), (G^1 v, \dots, G^n v)\right)\right],
    \label{eq:opt_1-to-n}
\end{equation}
where $v$ is drawn from certain probability distribution $\chi$ over the input space, and the cost function $C: \mathcal{V}^n \times \mathcal{V}^n \to \mathbb{R}$ is set to the relative mean square error: $C(a, b) = \|a - b\|_2 / \|b\|_2$.

In this work, we propose learning the extended $n$-step operator $\bar{G}$ using a neural network parameterization $\bar{\mathcal{G}}_{\bar{\theta}}$. The learning task is framed as the following minimization problem:
\begin{equation}
    \min_{ \bar{\theta} \in \bar{\Theta}} \mathbb{E}_{v \sim \chi} \left[C\left(\bar{\mathcal{G}}_{\bar{\theta}} v, (G^1 v, \dots, G^n v)\right)\right].
    \label{eq:opt_1-to-n2}
\end{equation}
After training, the recurrent application of the network $\bar{\mathcal{G}}_{\bar{\theta}}$ can be used to roll out predicted solutions of arbitrary length by iteratively updating the input function with the last element (i.e. the $n$-th element) of the previous output.

\section{Operator Learning Methods\label{sec:networks}}

In this section, we describe two baseline operator learning methods—Fourier Neural Operator (FNO) and Convolutional Neural Network (CNN)—for learning the single-stepping solution advancement operator $G$. We then introduce two extensions of these methods, the Koopman theory-inspired Fourier Neural Operator (kFNO) and Koopman theory-inspired Convolutional Neural Network (kCNN), which are designed to learn the multi-stepping solution advancement operator $\bar{G}$. 

\subsection{Baseline Fourier Neural Operator (FNO)}
The Fourier Neural Operator (FNO), introduced in \cite{FNO2020}, learns infinite-dimensional operators by parameterizing integral kernel operators in Fourier space. Applied to single-step solution advancement, FNO models the mapping $\disp(x,t) \mapsto \disp(x,t_1)$ using a sequence of transformations $\mathcal{P} \circ H \circ \mathcal{L}$: a lifting map $\mathcal{L}$, a hidden map $H$, and a projection map $\mathcal{P}$.

The lifting map $\mathcal{L}: \mathcal{V} \to \mathcal{V}_* ; \disp(x,t) \mapsto \z_0(x)$ elevates the input function into a higher-dimensional space $\mathcal{V}_* := \mathcal{V}_*(\mathbb{R}^d; \mathbb{R}^{d_\varepsilon})$, with $d_\varepsilon > d_v$. 

Next, the hidden map $H: \mathcal{V}_* \to \mathcal{V}_*$ performs a series of updates $\z_0 \mapsto \z_1 \mapsto \cdots \mapsto \z_L$, with each update $\z_l \mapsto \z_{l+1}$ implemented by a Fourier layer as follows:
\begin{equation}
\varepsilon_{l+1} = \alpha \varepsilon_{l} + \sigma\left( W_{l} \z_l + b_l +  \mathcal{F}^{-1} \{ \mathcal{R}_l( \mathcal{F} \{ \z_l \} ) \} \right),
\label{eq:FourierLayer}
\end{equation}
where $W_l \in \mathbb{R}^{d_\varepsilon \times d_\varepsilon} $ and $b_l \in \mathbb{R}^{d_\varepsilon}$ are learnable weights and biases, $\sigma$ is a nonlinear activation, and $\mathcal{F}$ and $\mathcal{F}^{-1}$ represent the Fourier transform and its inverse. The parameter $\alpha \in [0,1]$ controls the skip connection, where $\alpha=0$ disables it and $\alpha=1$ enables a full skip connection.
The operator $\mathcal{R}_l$ transforms the Fourier modes as
\begin{eqnarray}
 \mathcal{R}_l (\F\{\z\} )_{\kappa,i} =   
     \sum_{j=1}^{d_\z}  (R_l)_{\kappa,i,j} \F \{\z\}_{\kappa,j} , \nonumber \\
\kappa=0,1,\ldots,\kappa^{\text{max}}, \quad i=1,\ldots,d_\z
\end{eqnarray}
where $R_l \in  \mathbb{C}^{\kappa^{\text{max}} \times d_\z \times d_\z} $ is a learnable tensor. 
Finally, the projection map $\mathcal{P}: \mathcal{V}_* \to \mathcal{V} ; \z_L(x) \mapsto \disp(x,t_1)$ maps the high-dimensional hidden state back to the low-dimensional function $\disp(x,t_1)$.
On above, both $\mathcal{L}$ and $\mathcal{P}$ can be implemented as simple multilayer perceptrons (MLPs).

\subsection{Koopman Theory-Inspired Fourier Neural Operator (kFNO)}

\begin{figure*}
	\centerline{
		\includegraphics[width=0.9\linewidth]{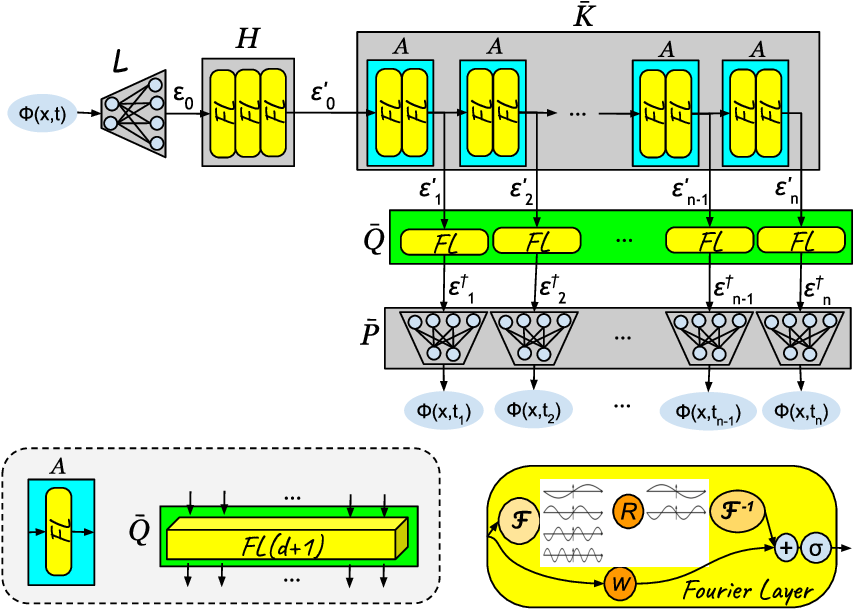}
	}
	\caption{
		\label{fig:kFNO}
		The Koopman theory-inspired Fourier Neural Operator model (kFNO). Note, replacing $A$ (or $\bar{Q}$) by the corresponding one in the left bottom box yields kFNO* (or kFNO$^\dagger$).
	}
\end{figure*}
The Koopman-inspired FNO (kFNO) extends the FNO architecture to learn a multi-step solution advancement operator $\bar{G}$. This is achieved by introducing a sequence of maps: $\bar{\mathcal{G}}_{\bar{\theta}} = \bar{\mathcal{P}} \circ \bar{Q} \circ \bar{K} \circ H \circ \mathcal{L}$ (see Figure~\ref{fig:kFNO}).

Similar to the baseline FNO, the lifting map $\mathcal{L}: \mathcal{V} \to \mathcal{V}_* ; \disp(x,t) \mapsto \z_0(x)$ and the hidden map $H: \mathcal{V}_* \to \mathcal{V}_* ; \z_0(x) \mapsto \z'_0(x)$ elevate the input into a higher-dimensional space and update it sequentially. 

The Koopman-inspired step involves the map $\bar{K}: \mathcal{V}_* \to \mathcal{V}_*^n ; \z'_0(x) \mapsto (\z'_1(x), \z'_2(x), \ldots, \z'_n(x))$, where $\z'_j(x) = A^j \z'_0(x)$, with $A: \mathcal{V}_* \to \mathcal{V}_*$ being an advancement operator. This allows kFNO to handle multiple time steps by advancing the state iteratively. 

The map $\bar{Q}: \mathcal{V}_*^n \to \mathcal{V}_*^n ; (\z'_1(x), \ldots, \z'_n(x)) \mapsto (\z^\dagger_1(x), \ldots, \z^\dagger_n(x))$ may be implemented by two variants: either using stacked Fourier layers that update $\z'_i(x) \mapsto \z^\dagger_i(x)$ for each $i = 1, \ldots, n$, or using one-higher-dimensional Fourier layers treating all inputs as a function having one extra domain dimension, i.e. $\z^\dagger: \mathbb{R}^{d+1} \to \mathbb{R}^{d_\z}$ .

The projection map $\bar{\mathcal{P}} $ then converts the multi-step hidden states  back to the solution space, yielding the output sequence $(\disp(x,t_1), \ldots, \disp(x,t_n))$. This is implemented as $n$ updates of $e^\dagger_i(x) \mapsto \disp(x,t_i)$ for $i=1,\ldots,n$, using the same baseline map $\mathcal{P}$.

In this framework, the operator $A$ serves as a time-advancement operator in a high-dimensional space, akin to the Koopman operator. It can be implemented either as a linear Fourier layer (without nonlinear activation) or a more expressive nonlinear operator with two stacked Fourier layers.

\subsection{Baseline Convolutional Neural Network (CNN)}

The Convolutional Neural Network (CNN) provides an alternative approach for learning the operator $G$ as a mapping from discretized solutions to their next time step. The architecture follows an encoder-decoder structure, similar to that of U-Net \cite{UNet} and ConvPDE \cite{ConvPDE}.

Let $e_0$ denote the input function $\disp(x_j,t)$ represented on an $x$-mesh and assume $p$ number of points along each spatial dimension.
The encoder block follows an iterative update procedure: $e_{l} \mapsto e_{l+1}$ over the level sequence $l = 0, 1, \dots, L-1$. Denote the last encoder output as $e'_L = e_L$, and apply a subsequent decoding procedure $(e'_{l+1}, e_l) \mapsto e'_l$ by reversing the levels.

Here, $e_l, e'_l \in \mathbb{R}^{c_l \times s_l}$ represent two data sequences, each with $c_l$ channels and size of $s_l$.
The data size starts from $s_0 = s_1= p^d$ and reduces to $s_l = (p/2^{l-1})^d$ as the level increases $l>1$.
 The encoder update $e_l \mapsto e_{l+1}$ is implemented using stacked convolutional layers (with a filter size of 3, stride 1, periodic padding, and ReLU activation). Some layers are replaced by Inception layers for improved performance. Additionally, a size-2 max-pooling layer is prepended to halve the size along each dimension for $l \geq 1$. 

The decoder update $(e'_{l+1}, e_l) \mapsto e'_l$ involves concatenating $e'_{l+1}$ (upsampled with size 2) with $e_l$ along the channel dimension. The final output is obtained as $\disp(x_j, t_1) = e'_1$.

\subsection{Koopman Theory-Inspired Convolutional Neural Network (kCNN)}

\begin{figure*}
    \centerline{
        \includegraphics[width=1\linewidth]{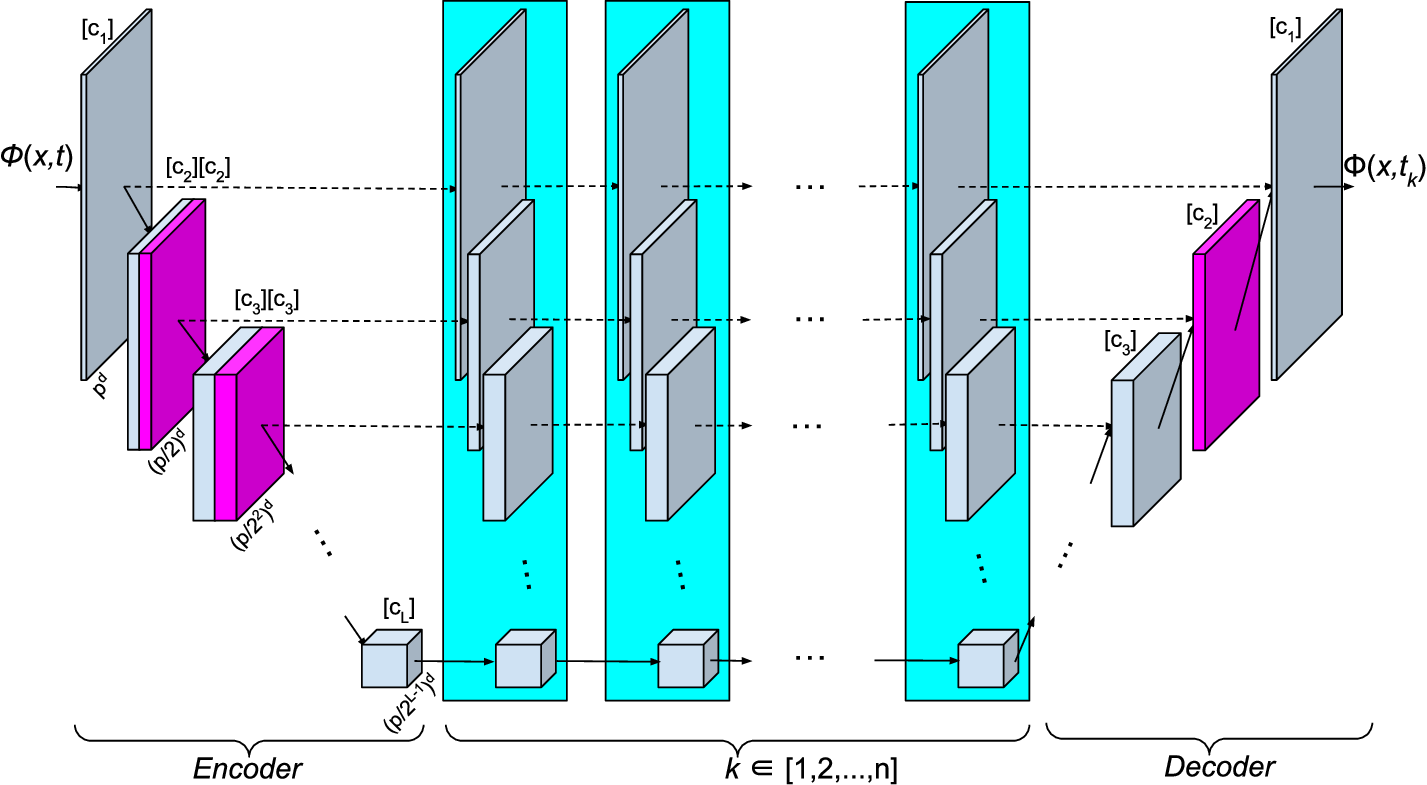}
    }
    \caption{
    \label{fig:CNN}
    The Koopman-inspired CNN model for learning the extended solution time advancement map $\disp(x,t) \mapsto (\disp(x,t_1), \dots, \disp(x,t_n))$. Standard convolutional layers are represented by gray blocks, which process $d$-dimensional input with a total pixel count of $(p/2^l)^d$ at different levels $l = 1, \dots, L$. The magenta block denotes convolutional layer variants, such as Inception layers. The output data channels $c_l$ for each convolutional layer are indicated within brackets.
    }
\end{figure*}

Similar to kFNO, the Koopman-inspired CNN (kCNN) learns a multi-step solution advancement operator $\bar{G}: \disp(x,t) \mapsto (\disp(x,t_1), \dots, \disp(x,t_n))$ by extending the baseline CNN architecture (see Figure~\ref{fig:CNN}).

The encoder in kCNN maps the input $\disp(x_j, t)$ to hidden representations $(e_1, \dots, e_L)$. Instead of immediately decoding these to predict the next step, kCNN introduces a hidden Koopman-like map $\mathcal{A}: (e_0, \dots, e_L) \mapsto (e^*_0, \dots, e^*_L)$, where $e^*_l$ represents the advanced states. 
 The map $\mathcal{A}^k$ allows the model to propagate the hidden states through multiple time steps. These updated representations are then fed into the decoder to generate the solution at subsequent time steps $(\disp(x_j,t_1), \dots, \disp(x_j,t_n))$.

The hidden map $\mathcal{A}$ is implemented by convolutional layers applied at each level of the hidden representations. This structure enables kCNN to learn time-evolving dynamics over multiple steps while leveraging the convolutional architecture’s strengths in handling spatial patterns.

\section{Results and Discussion \label{sec:results} }

In this section, we present our findings on the application of four operator learning models - kFNO, kCNN, baseline FNO, and baseline CNN - to predict flame evolution driven by two fundamental instabilities: Darrieus-Landau (DL) \cite{DARRIEUS1938UNPB,landau1988theory} and Diffusive-Thermal (DT) \cite{zeldovich1944selected,sivashinsky1977diffusional}. These instabilities are modeled using the $d$-dimensional Sivashinsky equation \cite{SivaEq}.

We generate training and validation datasets through high-order numerical simulations of the Sivashinsky equation for both one-dimensional ($d=1$) and two-dimensional ($d=2$) cases. These simulations use parameter conditions relevant to modeling the development of DL and DT fronts. We then train the models to learn the solution advancement operator for both types of fronts.

The performance of each model is evaluated based on its ability to predict both short-term and long-term flame evolution across different parameter regimes. Our assessment focuses on how well the models capture the characteristic behaviors of each instability type.

In the following subsections, we first describe the governing equations that form the basis of our study. We then discuss the training datasets and present a detailed analysis of each model's performance in capturing the flame evolution dynamics. Finally, we compare the models' effectiveness in learning both one-dimensional and two-dimensional flame instability development.

\subsection{Governing Equations of Flame Instability}

The development of a statistically planar flame front under intrinsic instabilities is modeled using a displacement function $\hat{\phi}(\hat{x}, \hat{t}): \mathbb{R}^d \times \mathbb{R} \to \mathbb{R}$, which describes the flame's stream-wise coordinate. The evolution of this function is governed by the $d$-dimensional Sivashinsky equation \cite{SivaEq}:

\begin{equation} 
\hat{\phi}_{\hat{t}} + \frac{1}{2} \left( \nabla \hat{\phi} \right)^2 
=
 -4(1 + \text{Le}^*)^2  \nabla^4 \hat{\phi} - \text{Le}^* \nabla^2 \hat{\phi} + (1 - \Omega) \Gamma(\hat{\phi})
\label{eq:MKS_org}
\end{equation}

Here, $\Gamma: \hat{\phi} \mapsto \mathcal{F}^{-1} \left( |\kappa| \mathcal{F}_\kappa (\hat{\phi}) \right)$ is a linear singular non-local operator defined using the spatial Fourier transform $\mathcal{F}_\kappa(\hat{\phi})$ and its inverse $\mathcal{F}^{-1}$, with $\kappa$ representing the wavevector. The operators $\nabla$, $\nabla^2$, and $\nabla^4$ denote the spatial gradient, Laplacian, and biharmonic operators, respectively. The parameter $\Omega$ represents the density ratio between burned products and fresh reactants, while $\text{Le}^*$ is determined by the Lewis number of the deficient reactant and a critical Lewis number.

Through appropriate rescaling of time $t \sim \hat{t}$, space $x \sim \hat{x}$, and the displacement function $\hat{\phi}(\hat{t},\hat{x}) \sim \phi(t,x)$, along with suitable parameter choices, the Sivashinsky equation can be reduced to two distinct limiting cases \cite{YHN2024}:

\begin{enumerate}
    \item The Michelson-Sivashinsky (MS) equation \cite{michelson1977nonlinear}, representing pure DL instability:

    \begin{equation}
    \frac{1}{\tau} \phi_t + \frac{1}{2 \beta^2} ( \nabla \phi )^2 = \frac{1}{\beta^2} 
    \nabla^2 \phi + \frac{1}{\beta} \Gamma(\phi)
    \label{eq:MS}
    \end{equation}

    \item The Kuramoto-Sivashinsky (KS) equation \cite{kuramoto1978diffusion}, describing pure DT instability:

    \begin{equation}
    \phi_t + \frac{1}{2 \beta^2} ( \nabla \phi )^2 = -\frac{1}{\beta^2} \nabla^2 \phi - \frac{1}{\beta^4} \nabla^4 \phi
    \label{eq:KS}
    \end{equation}
\end{enumerate}

In both equations, $\phi(x, t)$ represents the rescaled front solution within a periodic domain $x \in \mathcal{D} := (-\pi, \pi]^d \subset \mathbb{R}^d$. The parameter $\beta$ corresponds to the largest unstable wave number according to linear stability analysis, and $\tau = \beta / 10$ adjusts the timescales of the DL instability.

Before proceeding further, it may be worthwhile to note a few known results. The KS equation \eqref{eq:KS} is often utilized as a benchmark example for PDE learning studies and is renowned for exhibiting chaotic solutions at large $\beta$. On the other hand, the MS equation \eqref{eq:MS}, although less familiar outside the flame instability community, can be precisely solved using a pole-decomposition technique \cite{Thual_Frisch_Henon_poledecomp}, transforming it into a set of ODEs with finite freedom. Moreover, at small $\beta$, the MS equation admits a stable solution in the form of a giant cusp front. However, at large $\beta$, the equation becomes susceptible to noise, resulting in unrest solutions characterized by persistent small wrinkles atop a giant cusp. Details about known theory can be found in references \cite{Vaynblat_matalon_polestability1, Vaynblat_matalon_polestability2, Olami_noise, denet2006stationary, Kupervasser_pole_book, Karlin2002cellular, Creta2020propagation, CRETA2011INST,rasool2021effect,YBB15PRE}.

\subsection{Training Dataset \label{sec:training_data}}

We employ a pseudo-spectral approach \cite{kassam2005fourth} to numerically solve Equations \eqref{eq:MS} and \eqref{eq:KS}. The linear terms are transformed into Fourier space, and the solution time advancement due to these terms is directly handled using the integrating factor method. The full solution advancement, accounting for the remaining nonlinear terms, is performed using a Runge-Kutta (4,5) time integration scheme.
 
For the one-dimensional ($d=1$) cases, all solutions are computed on a uniformly spaced mesh consisting of 256 points. We consider a total of four cases: two for the DL instability and two for the DT instability, with $\beta = 10$ and $\beta = 40$. For each case, we generate 250 sequences of short-duration solutions, each spanning a time duration of $0 \leq t \leq 75$ and containing 500 consecutive solutions separated by a time interval of $\Delta_t = 0.15$.
Note that with a larger time step, the 1D pseudo-spectral solver tends to crash when simulating the DL case at $beta=40$ due to the aforementioned solution sensitivity to uncontrolled noise.
Additionally we generate a single sequence of long-duration solutions, covering a time duration of $0 \leq t \leq 18,750$ and comprising 125,000 consecutive solutions outputted at the same interval $\Delta_t$.
 Each short solution sequence starts from random initial conditions $\phi_0(x)$ sampled from a uniform distribution over the range $[0, 0.03]$ in physical space. A validation dataset is similarly created for all four cases, containing 10\% of the data present in the training dataset.

 For the two-dimensional ($d=2$) cases, solutions are computed on a uniform mesh of size 128$^2$. Two cases are considered: one for the DL instability and one for the DT instability, both with $\beta = 15$. For each case, 45 sequences of medium-duration solutions are generated, each spanning the time interval $0 \leq t \leq 435.12$ and consisting of 5880 consecutive solutions with a time step of $\Delta_t' = 0.074$. This time step is again determined to ensure stable simulation by the spectral solver on a coarser 2D mesh.  Each solution sequence begins from a random initial condition $\phi'_0$ that contains only low wave-number fluctuations. Specifically, wave content beyond eight wave numbers is truncated to zero, i.e., $\mathcal{F}_{\kappa}\{\phi'_0\} = 0 , \forall  |\kappa| > 8$.  Note that in \cite{Yu2023}, two strategies for randomizing initial conditions in wave space and physical space were used, which yielded similar model learning performance.
 A validation dataset, containing one ninth of the data in the training dataset, is created for both two cases.

\subsection{Analysis of One-Dimensional Results}

\begin{figure*}
	\centerline{
	\includegraphics[width=1\linewidth]{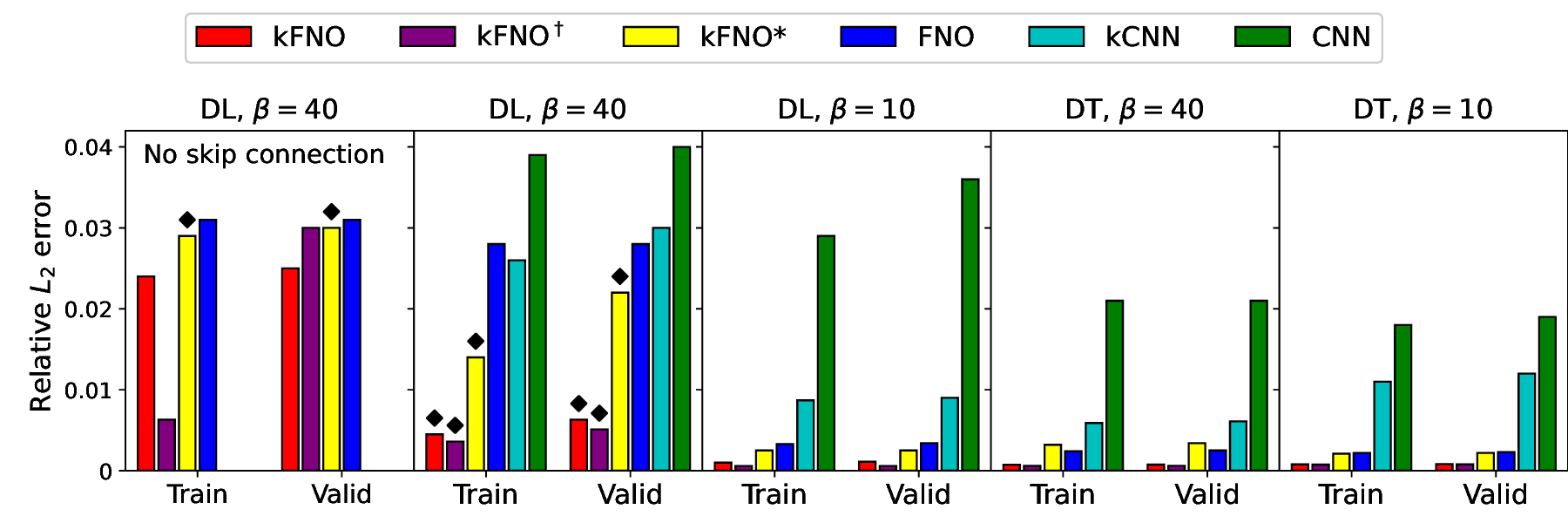}
	}
	\caption
	{
		\label{fig:errTable}
		Short-term relative $L_2$ training/validation errors for all models in learning 1D instabilities. 
		The first subfigure correpond to FNO-based models without skip connections in the Fourier layers. 
		The diamiond symbol indicate models prone to divergence in long-term recurrent predictions. Detailed values can be read from Table \ref{Table1} in the appendix.		
	}
\end{figure*}

We trained the solution advancement operators $G : \phi(x;t) \mapsto \phi(x;t_1)$ and $\bar{G} : \phi(x;t) \mapsto ( \phi(x;t_1), \ldots, \phi(x;t_n) )$ with $n=20$ using the training datasets for DL and DT instabilities described in the section \ref{sec:training_data}.

The baseline FNO and CNN models are compared against their Koopman-inspired counterparts: kFNO and kCNN. Implementation details for the CNN and kCNN models are provided in Appendix \ref{app:nn_detail}. For the FNO-based models, the default kFNO architecture follows the design shown in Figure \ref{fig:kFNO}, where maps $H$, $A$, and $\bar{Q}$ are implemented using three, two, and one Fourier layer, respectively, with skip connections in all Fourier layers (i.e. let $\alpha=1$ in Eq. \eqref{eq:FourierLayer}) .

To ensure a fair comparison, the baseline FNO model is recovered from the kFNO architecture by setting $n=1$  and merging the $\bar{K}$ and $\bar{Q}$ maps into the hidden map $H$. This approach guarantees that the baseline FNO and kFNO models have equivalent expressive power, allowing for a direct and unbiased comparison of their performance.

Additionally, two variants of kFNO are explored: kFNO* and kFNO$^\dagger$. In kFNO*, the hidden map $A$ is replaced by a single linear Fourier layer (without the last nonlinear activation). In kFNO$^\dagger$, the map $\bar{Q}$ is replaced by a high-dimensional Fourier layer that treats the $n$ $x$-dependent inputs as a single function. These modifications are highlighted in the bottom left corner of Figure \ref{fig:kFNO}.

The default kFNO model, which includes skip connections in all Fourier layers, significantly reduces training and validation errors compared to models without skip connections (i.e. letting $\alpha=0$ in Eq. \ref{eq:FourierLayer}). For the two DT cases as well as the DT case with small $\beta=10$, the trained kFNO models accurately predict long-term solutions with correct statistics. However, for the DL case of large $\beta=40$ , the solution is sensitive to numerical round-off errors, causing the trained models to diverge after a few hundred recurrent prediction steps, making them unsuitable for long-term predictions.

Interestingly, for this challenging case, removing the skip connections in the Fourier layers enabled models to predict long-term solutions with accurate statistics. However, this improvement in long-term performance comes at the cost of short-term accuracy, as these models exhibit higher training and validation errors compared to those with skip connections. In the following, when we refer to kFNO for the challenging DL case ($\beta=40$) without additional clarification, we are referring to the model without skip connections in the Fourier layers. For all other cases, kFNO refers to the default model with skip connections.

\begin{figure*}
	\centerline{
	\includegraphics[width=1\linewidth]{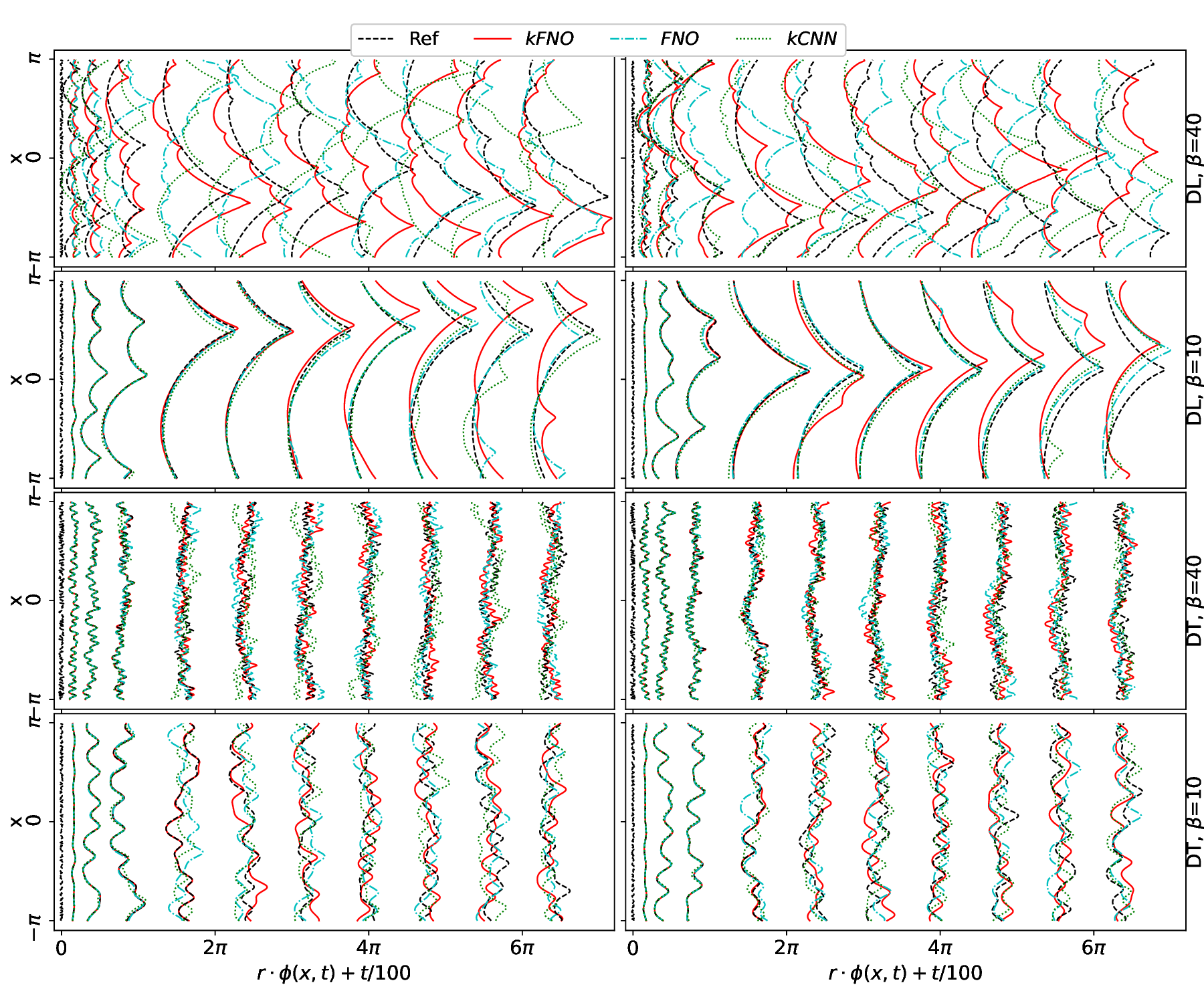}
	}
	\caption
	{
		\label{fig:disp_1D}
Long-term solutions of 1D flame front displacement $\disp(x,t)$ for four instability cases (DL/DT at $\beta = 40$ and $\beta = 10$, shown in rows). Black reference solutions from the spectral solver are compared to model predictions by kFNO (red), FNO (cyan), and kCNN (green). The left and right columns represent two randomly initialized solution sequences, each displaying eleven snapshots of $\disp(x,t)$ at $t/\Delta_t$= 0, 50, 125, 250, 500, 750, 1000, 1250, 1500, 1750, and 2000. For better visualization, the fronts are rescaled by a factor $r$ and shifted by $t/100$ to avoid overlap. The scaling factor is $r = 1$, except for the third-row DT case with $\beta = 10$, where $r = 2$.
	}
\end{figure*}

\begin{figure*}
	\centerline{
		\includegraphics[width=1\linewidth]{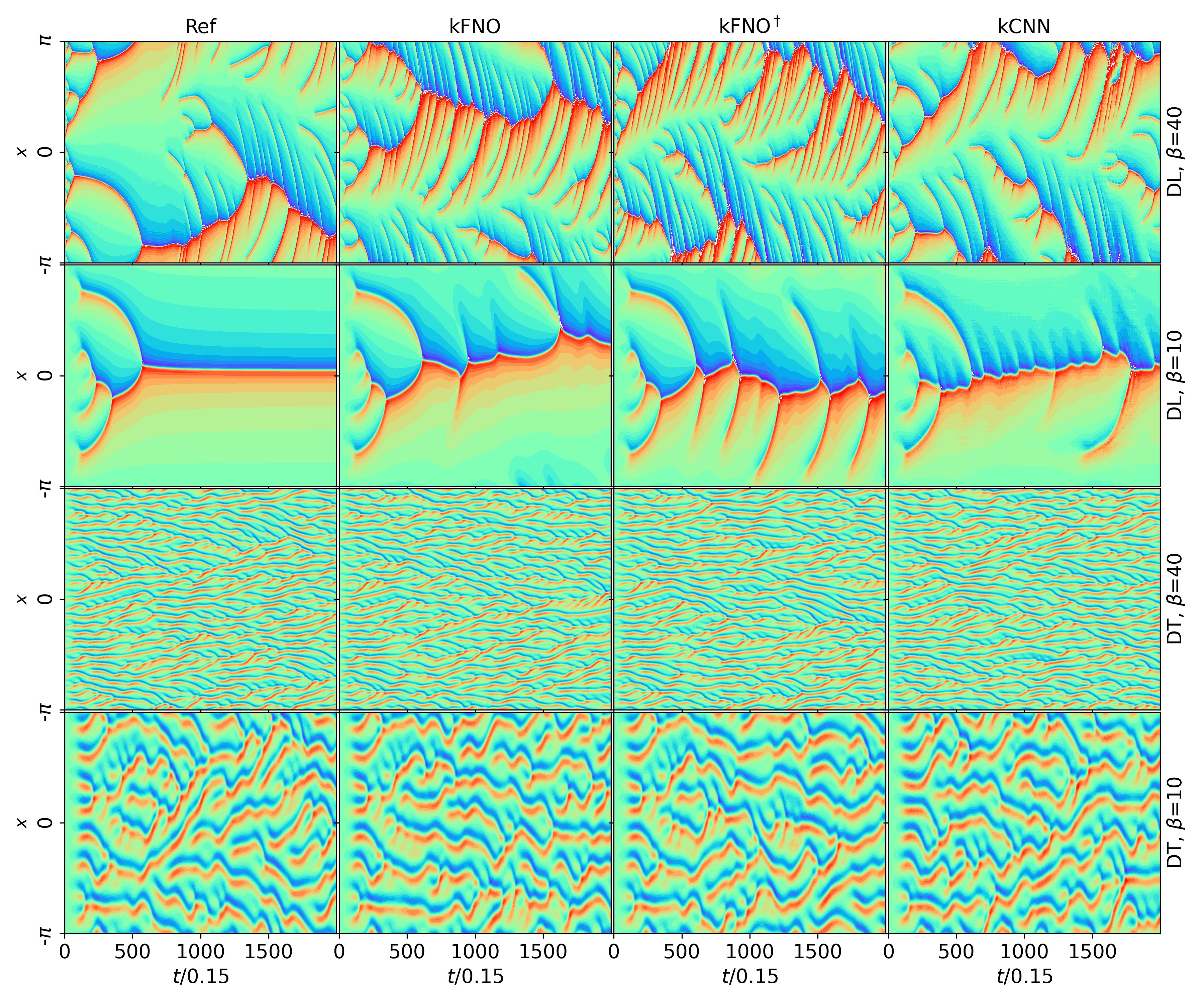}	   
}
	\caption{
		\label{fig:uSlope_1D}
		Comparison of 1D front slope $\disp_x(x,t)$ for four instability cases (DL/DT at $\beta = 40$ and $\beta = 10$, shown in different rows), between a reference solution (first column) and predictions from kFNO, kFNO$^\dagger$, and kCNN (last three columns). The rainbow color scale represents values from negative to positive.
		}
\end{figure*}

\begin{figure*}
	\centerline{
		\includegraphics[width=1\linewidth]{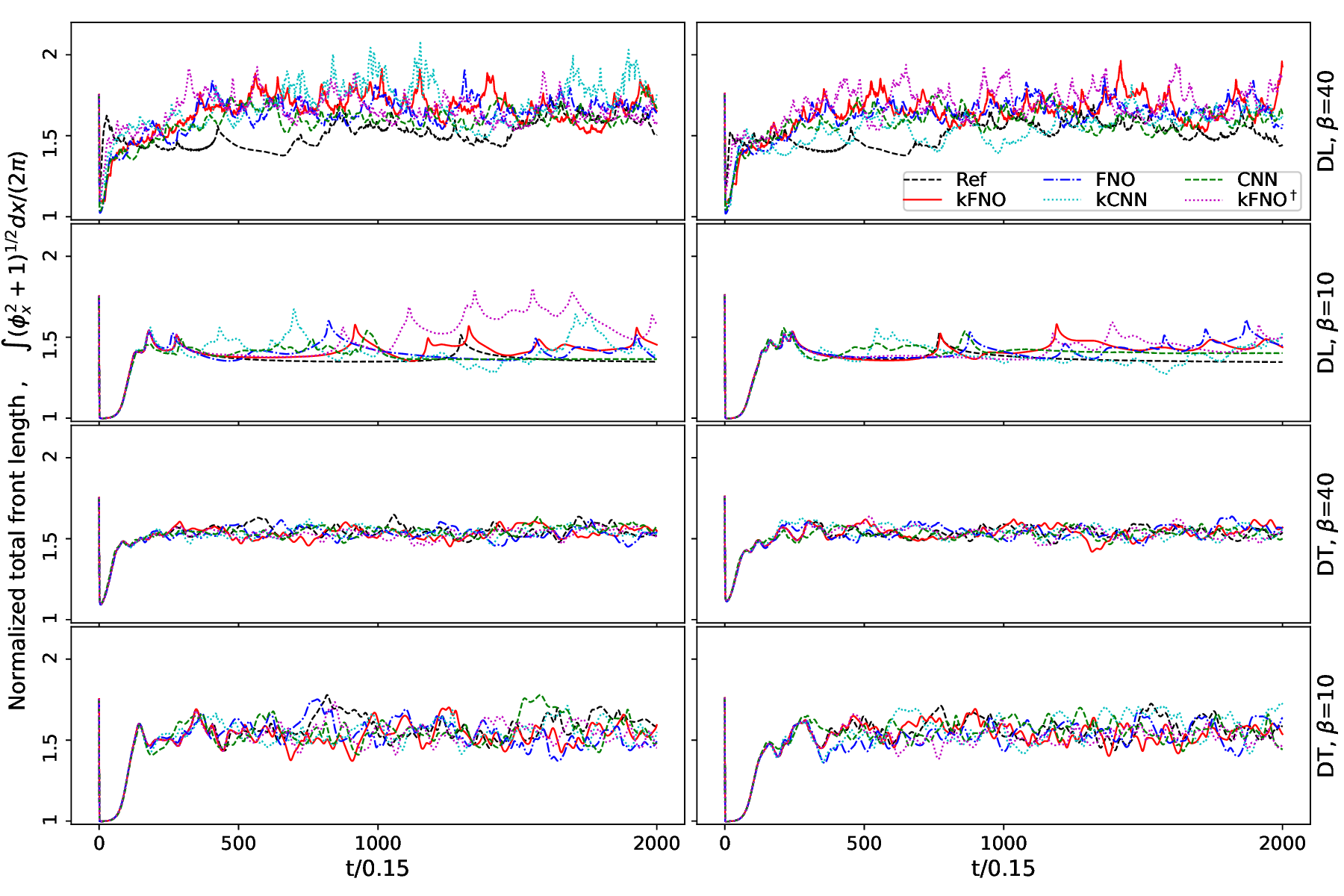}
   }
	\caption{
		\label{fig:len_1D}
Normalized total front length comparison at four 1D instability cases (DL/DT at $\beta = 40$ and $\beta = 10$, shown in rows).
Two instances of randomly initialized cases are shown in the left and right columns. The black curve represents reference solutions obtained from spectral solver. Predictions by kFNO, FNO, kCNN and CNN are shown in red, blue, cyan and green respectively.
	}
\end{figure*}

\begin{figure*}
	\centerline{
		\includegraphics[width=1\linewidth]{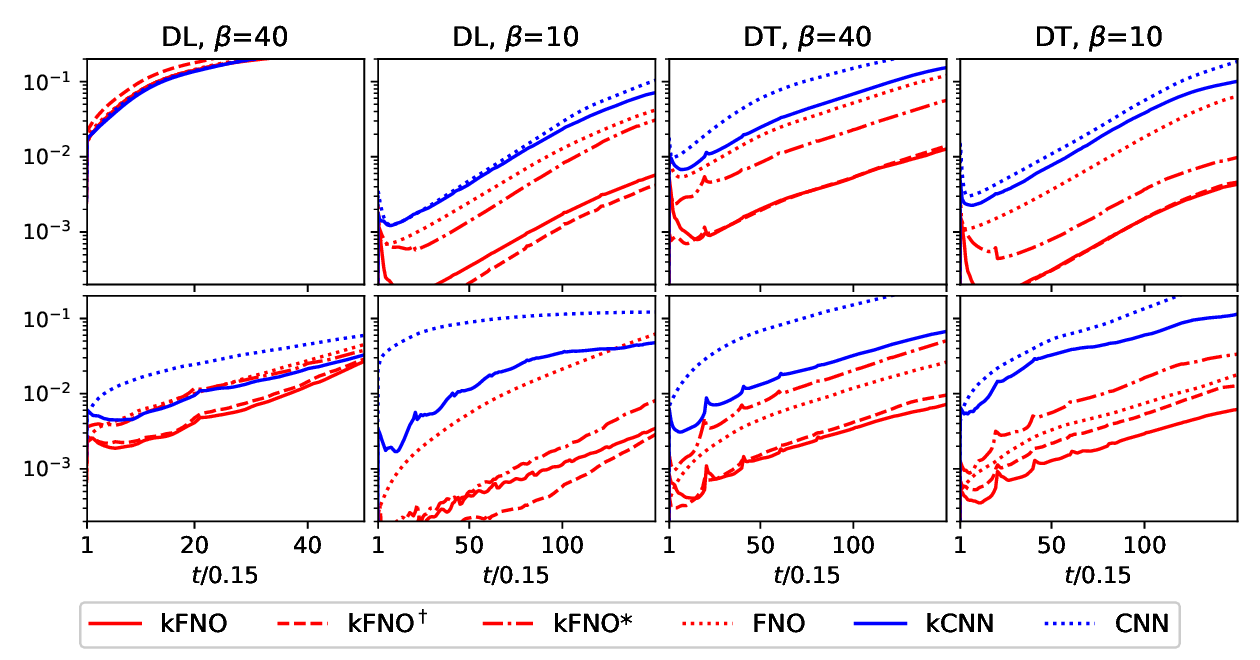}
   }
	\caption{
		\label{fig:ave_err_1D}
		Time evolution of spatially-averaged relative $L_2$ error  between reference 1D solutions and predictions by kFNO,  kFNO$^\dagger$, kFNO*, FNO, kCNN, and CNN (see legends) for four cases (DL/DT at $\beta = 40$ and $\beta = 10$ in columns).
		The error in the first row is averaged over 10 cases starting with random initial conditions, while the second row shows the error averaged over 10 cases starting from long-term reference solutions evolved after 2000 time steps with a step size of 0.15.
		}
\end{figure*}

\begin{figure*}
	\centerline{
		\includegraphics[width=1\linewidth]{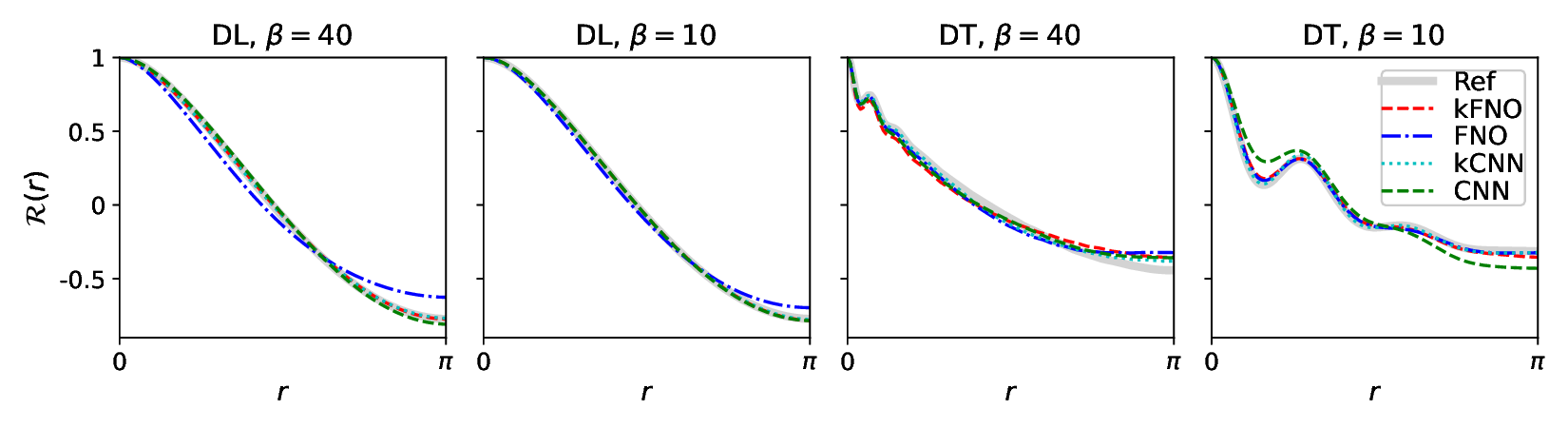}
   }
	\caption{
		\label{fig:corr_1D}
		 Comparison of the auto-correlation function $\mathcal{R}(r)$, Eq. \eqref{eq:corr}, which measures 1D long-term solution statistics, between reference solutions and predictions by kFNO, FNO, kCNN, and CNN for four cases (DL/DT at $\beta = 40$ and $\beta = 10$, shown in columns).
		 The expectation in Eq. \eqref{eq:corr} is estimated by averaging over 10 random cases for the time duration $500 < t/0.15 < 4000$.		 
		}
\end{figure*}

Fig. \ref{fig:errTable} (or table \ref{Table1} in the appendix ) summarizes the relative L2 training/validation errors for all models. Additionally, Figure \ref{fig:disp_1D} compares reference front displacement solutions $\disp(x,t)$ with predictions from kFNO, FNO, and kCNN for four cases (see \cite{Yu2023} for similar plot by CNN). Figure \ref{fig:uSlope_1D} compares the reference front slopes $\disp_x(x,t)$ to those predicted by kFNO, kFNO$^\dagger$ and kCNN (see \cite{Yu2023} for similar plot by FNO and CNN), while Figure \ref{fig:len_1D} compares the total front area predictions by kFNO, kFNO$^\dagger$, FNO, kCNN, and CNN. 
Figure \ref{fig:ave_err_1D} presents the temporal development of the spatially-averaged relative L2 errors between reference solutions and predictions by the four models of kFNO, FNO, kCNN, and CNN (the non-averaged errors are shown in Figure \ref{fig:err_1D} of Supplementary material).
The long-term statistical behavior of the predicted 1D solutions can be quantified using the auto-correlation function:
\begin{equation}
\mathcal{R} (r)   = \mathbb{E} \left(
 \int_\D \phi^*(x) \phi^*(x-r) dx 
 /\int_\D \phi^*(x) \phi^*(x) dx 
 \right).
\label{eq:corr}
\end{equation}
where $\phi^*(x)$ denotes a predicted solution after a sufficiently long time. Figure \ref{fig:corr_1D} compares the auto-correlation function obtained from the four models with the reference solution.

The following key observations can be made:

\textbf{Short-term accuracy}: The kFNO and kCNN models provide highly accurate short-term predictions for front solutions across all DL and DT cases. As shown in figure \ref{fig:errTable}, both models consistently achieve significantly lower training and validation errors compared to their baseline FNO and CNN counterparts. Notably, all four kFNO models, which incorporate skip connections in the Fourier layers, exhibit error rates that are three times lower than those of the corresponding baseline FNO models. Similarly, the kCNN error is three times lower than the CNN one for the two cases (DT at $\beta=40$) and (DL at $\beta=10$). This reduction is further highlighted by the error growth comparison between kFNO and FNO, and between kCNN and CNN, as depicted in Figure \ref{fig:ave_err_1D}.

\textbf{Variant performance}: The kFNO* model, which uses a single linear Fourier layer for the hidden map $A$, tends to have larger training and validation errors than the default kFNO, although it still outperforms the baseline FNO model. Except for the  DL case of large ($\beta=40$), kFNO* is capable of predicting long-term solutions with accurate statistics.

On the other hand, kFNO$^\dagger$, which uses a high-dimensional Fourier layer in the $\bar{Q}$ map, shows comparable training and validation errors to the default kFNO model. However, kFNO$^\dagger$ demonstrates a tendency to overfit in the most challenging DL flame case ($\beta=40$), where the training error is considerably lower than the validation error (fig. \ref{fig:errTable}). Nevertheless, this model can predict long-term solutions and reproduce accurate statistics for all flame cases (Figures \ref{fig:uSlope_1D} and \ref{fig:len_1D}) .

\textbf{FNO vs. CNN}: In general, FNO-based models outperform CNN-based models, a well-known result from previous studies. This holds true for the Koopman-inspired variants as well, with kFNO and kCNN inheriting the improved accuracy of FNO.

\textbf{Long-term performance}: Both kFNO and kCNN models successfully predict long-term solutions and reproduce correct statistics, as evidenced by the comparisons of predicted front displacement (Figure \ref{fig:disp_1D}), front slope (Figure \ref{fig:uSlope_1D}), total flame length (Figure \ref{fig:len_1D}) and auto-correlation function (Figure \ref{fig:corr_1D}). For the two challenging DL flame cases, both models capture the sensitivity of the solution to noise, particularly at $\beta=40$. However, all models slightly overestimate the noise effect in the DL cases of small $\beta=10$, where the analytical solution tends toward a steady state. This behavior has also been reported in previous studies for baseline models.

\textbf{Computation speedup}: Both kFNO and kCNN exhibit greater computational efficiency compared to their baseline counterparts when applied for long-term predictions. This speedup is achieved because the $n$-time-step solution advancement in the Koopman-inspired models operates in a high-dimensional space by repeatedly applying a simple operator (i.e., $A$ for kFNO and $\mathcal{A}$ for kCNN). In contrast, baseline models require additional computations at each time step, where the high-dimensional hidden state is projected back to the low-dimensional solution space and then lifted back to the hidden state, slowing down the prediction process.

\subsection{Analysis of Two-Dimensional Results}

\begin{table}
	\caption{ 
		Short-term relative $L_2$ training/validation errors for two models learning 2D instabilities governed by Eqs \eqref{eq:MS} and \eqref{eq:KS} with $\beta=15$.
		\label{Table2} }  
	\centerline{
		\begin{tabular}{ |c |c | c| }  
			\hline
		       & DL: Train L2/Valid. L2  &  DT: Train L2/Valid. L2  \\
			\hline
			2D-kFNO       & 0.0021/0.0015  &  0.0025/0.0025  \\ 
		   2D-FNO        & 0.0041/0.0030   &     0.0053/0.0053  \\ 
		  \hline
		\end{tabular}
	}
\end{table} 
 
We now extend our analysis to two-dimensional ($d=2$) flame front instabilities, focusing on the DL and DT cases at $\beta = 15$. 
The training datasets for these cases are generated as described in Section \ref{sec:training_data}. 
Due to the significant computational resource requirements of 2D model training, we only apply the FNO and kFNO models to learn the solution advancement operators 
$G : \phi(x;t) \mapsto \phi(x;t_1)$ and $\bar{G} : \phi(x;t) \mapsto (\phi(x;t_1), \ldots, \phi(x;t_n))$ with $n=20$, respectively. Note that the 2D time step of $0.074$ is half of the 1D time step. 
Except for minor details provided in Appendix \ref{app:nn_detail}, the 2D model architectures and hyperparameters are identical to those used for the 1D models.

\begin{figure*}
	\centerline{
		\includegraphics[width=1\linewidth]{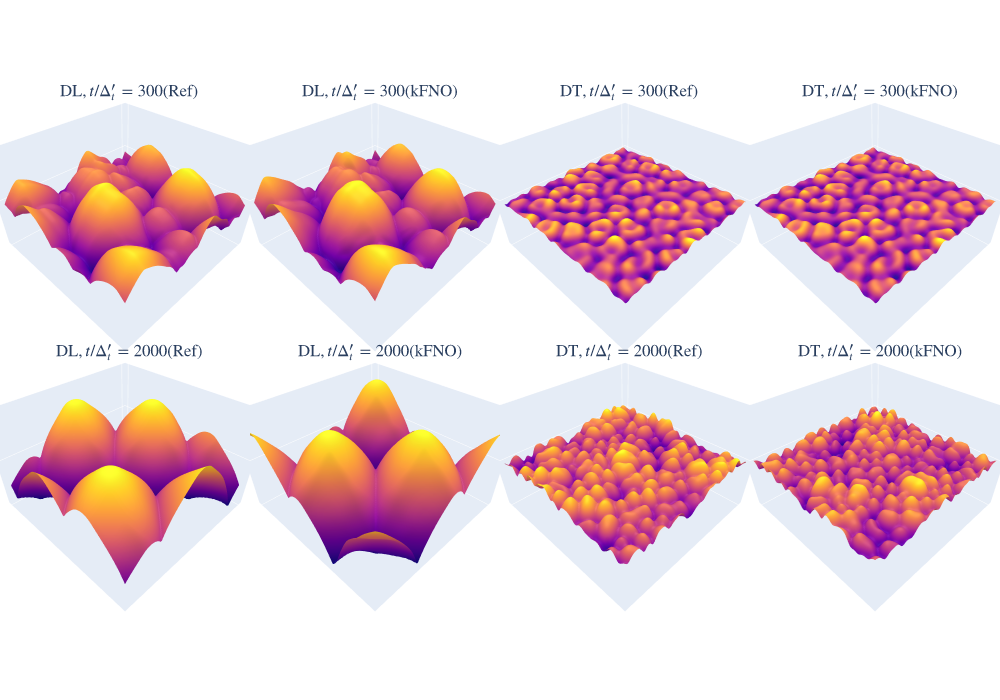} 
		}
	\caption{
		\label{fig:2DSnapshots}
		Comparison of two-dimensional snapshots of DL and DT fronts with $\beta = 15$, obtained from the reference spectral solver and kFNO predictions. The top row shows snapshots at $t/\Delta_t' = 300$, while the bottom row corresponds to $t/\Delta_t' = 2000$. Here, $\Delta_t' = 0.074$.
		}
\end{figure*}

\begin{figure*}
	\centerline{
		\includegraphics[width=1\linewidth]{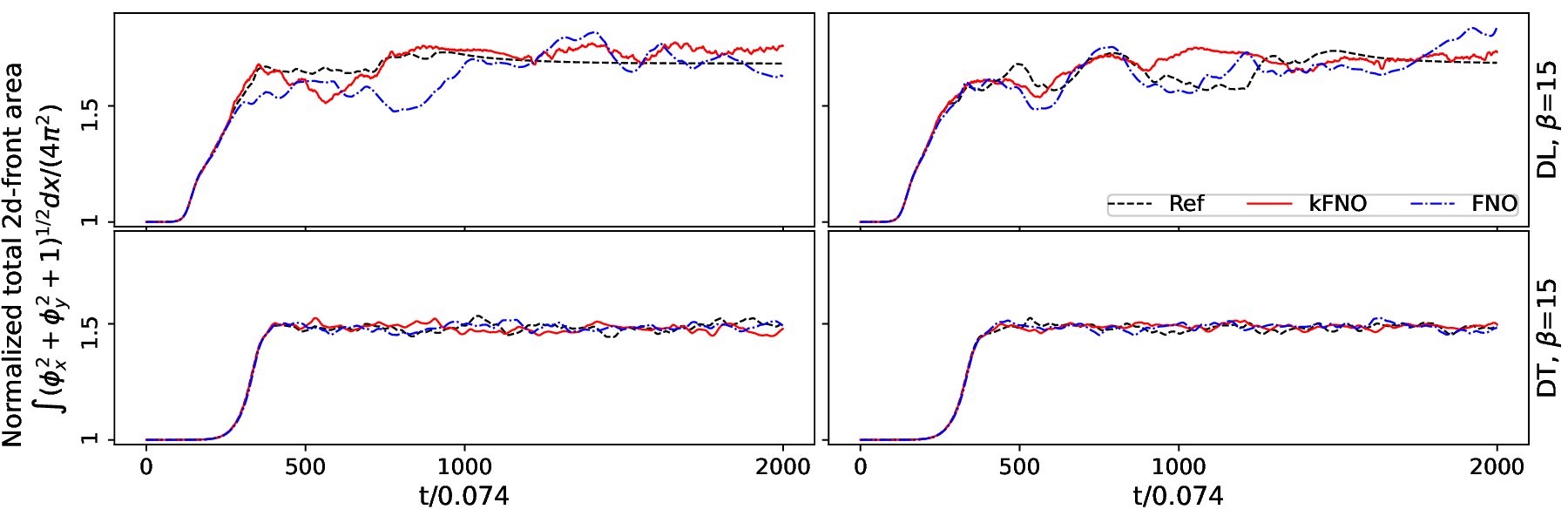} 
		}
	\caption{
		\label{fig:2Dlen}
		Comparison of two random instances (columns) of 2D front areas between the reference spectral solver (black), kFNO predictions (red), and baseline FNO predictions (blue) for four 2D instability cases (DL/DT at $\beta = 15$, shown by rows). Note: different from 1D case, the 2D initial conditions are randomly generated only allowing low wave-number fluctuations.
		}
\end{figure*}

\begin{figure*}
	\centerline{
		\includegraphics[width=0.65\linewidth]{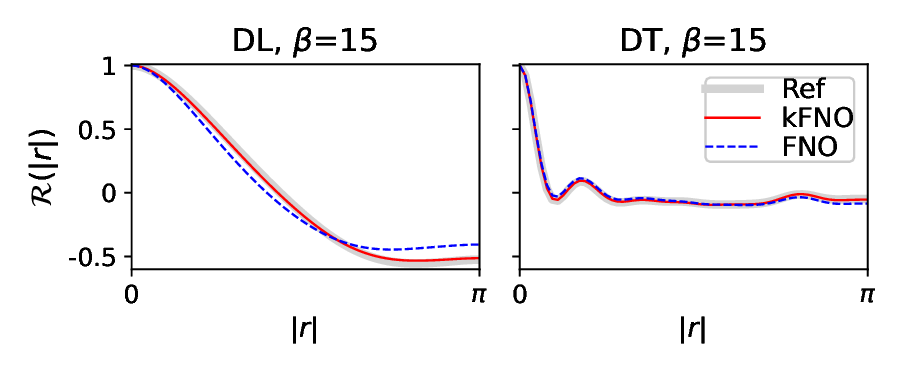}
   }
	\caption{
		\label{fig:corr_2D}
		 Comparison of the directional-independent auto-correlation function $\mathcal{R}(|r|)$ calculated for long term 2D solutions obtained by reference spectral solver and ones predicted by kFNO and FNO for DL/DT case at $\beta = 15$ (shown in columns).
		 The expectation in Eq.~\eqref{eq:corr} is estimated by averaging over five random run-instances over a time duration $500 < t / 0.074 < 4000$.
		}
\end{figure*}

\begin{figure*}
	\centerline{
		\includegraphics[width=0.65\linewidth]{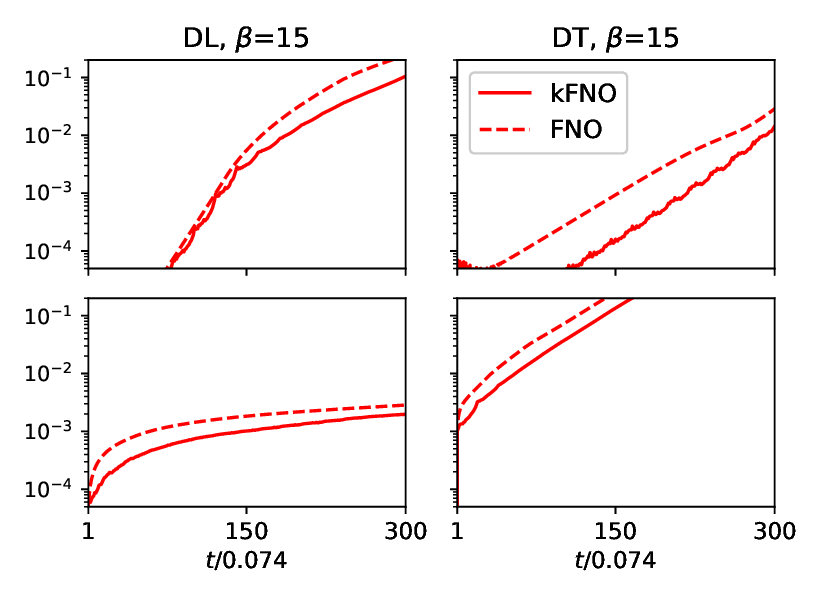} 
		}
	\caption{
		\label{fig:2Derr}
		Time evolution of spatially-averaged relative $L_2$ error between reference 2D solutions and kFNO/FNO predictions for DL/DT cases at $\beta = 15$. Top: errors averaged over 5 random runs with low-wavenumber initial conditions. Bottom: errors averaged over 5 random runs starting from long-term solutions after 2000 steps ($\Delta_t' = 0.074$).
		}
\end{figure*}

Table \ref{Table2} presents the relative $L_2$ training and validation errors for all models. 
Figures \ref{fig:2DSnapshots} compare typical snapshots of 2D DL/DT fronts obtained from the reference spectral solver with those predicted by kFNO. 
Figure \ref{fig:2Dlen} compares the evolution of the total front area from two randomly initialized reference 2D solutions with predictions by kFNO and FNO for both DL and DT cases. 
Figure \ref{fig:corr_2D} contrasts the direction-independent autocorrelation function $\mathcal{R}(|r|)$, computed from the long-term 2D solutions produced by the reference spectral solver and predictions by kFNO and FNO for the two cases. 
Figure \ref{fig:2Derr} depicts the temporal development of the spatially-averaged relative $L_2$ errors between the reference 2D solutions and predictions by kFNO and FNO for both cases. A plot of non-averaged error growth is provided in Figure \ref{fig:2Derr} of the Supplementary Material.

The following observations can be made regarding the 2D results:

The kFNO model effectively predicts long-term 2D solutions for both DL and DT cases at $\beta = 15$, accurately capturing the characteristic behaviors of each instability type, as shown in Figure \ref{fig:2DSnapshots}. The model reliably reproduces the evolution of the total front area (Figure \ref{fig:2Dlen}) and the autocorrelation function (Figure \ref{fig:corr_2D}). Additionally, it maintains low relative $L_2$ errors throughout the prediction process, as depicted in Figure \ref{fig:2Derr}.

Consistent with the 1D findings, kFNO consistently outperforms the baseline FNO model, delivering superior performance for both short-term and long-term predictions. Specifically, kFNO achieves two times smaller training and validation errors for both cases (see Table \ref{Table2}) and produces more accurate autocorrelation statistics, as demonstrated by the red and blue lines in Figure \ref{fig:corr_2D}.
Moreover, as shown in the appendix, the training of kFNO is significantly faster than that of the baseline FNO method on the same 2D dataset.

\section{Summary and Conclusion \label{sec:conclusion} }

This work presents a novel approach for learning solution advancement operators of nonlinear partial differential equations (PDEs) by integrating Koopman operator theory into Fourier Neural Operator (FNO) and Convolutional Neural Network (CNN) architectures. These extensions, named Koopman-inspired FNO (kFNO) and Koopman-inspired CNN (kCNN), are designed to capture both short-term accuracy and long-term statistical fidelity in time evolution, addressing the challenges posed by chaotic dynamics.

The proposed methods aim to learn the time advancement operator in a high-dimensional latent space, enabling accurate multi-step predictions. We benchmark the performance of kFNO and kCNN against baseline FNO and CNN models using two nonlinear PDEs modeling flame instabilities—the Michelson-Sivashinsky (MS) equation for hydrodynamic instabilities and the Kuramoto-Sivashinsky (KS) equation for diffusive-thermal instabilities. Performance metrics include short-term prediction accuracy and long-term statistics in one- and two-dimensional configurations.

Our results demonstrate that the Koopman-inspired frameworks outperform baseline methods in short term prediction accuracy and long-term statistical fidelity. Specifically, the kFNO and kCNN models consistently achieve lower training and validation errors compared to their FNO and CNN counterparts, respectively. The Koopman-inspired models not only accurately predict long-term solutions but also capture the correct statistical behavior of chaotic dynamics, as evidenced by the auto-correlation function. Furthermore, the Koopman-inspired models exhibit superior computational efficiency, enabling faster long-term predictions compared to baseline models. These findings emphasize the potential of operator learning methods rooted in Koopman theory for advancing computational modeling of complex physical systems.

Future research directions include extending the proposed methods to more complex PDEs under different boundary conditions and exploring the integration of Koopman operator theory into other architectures such as DeepOnet\cite{DeepONet}. Additionally, the application of Koopman-inspired models to real-world problems in fluid dynamics, climate science, and other fields is an exciting avenue for future investigation \cite{herbert2024comparison,hodzic2017a,hodzic2018,YBL15,YU2017nonlinear,YNBL2019_CF,YL2019_Equation,YNCL2020_JFM}.

\begin{acknowledgments}
The authors gratefully acknowledge the financial support by the Swedish Research Council(VR-2019-05648) and the AI Lund initiative grant. The stay abroad of M.H. was supported by the Federal Ministry of Defence.
The computations were enabled by resources provided by the National Academic Infrastructure for Supercomputing in Sweden (NAISS), at ALVIS and Tetralith, partially funded by the Swedish Research Council through grant agreement no. 2022-06725.
\end{acknowledgments}

\section*{Data Availability Statement}
The code and data that support the findings of this study are openly available at www.github.com/RixinYu/ML\_paraFlame.

\section*{Declaration of Interest}
The authors have no conflicts to disclose.


\appendix
\section{ Model hyper-parameters and training details \label{app:nn_detail} }

All models are trained for 1000 epochs using the Adam optimizer with a learning rate of 0.0025 and a weight decay of 0.0001. A learning rate scheduler is applied with a step size of 100 and a decay factor (gamma) of 0.5. To ensure training stability, gradient clipping is employed with a maximum norm threshold of 30.

For 1D models, training is conducted with a batch size of 1000 on a single GPU (NVIDIA A40). All 1D training processes are completed within 4 hours. For 2D models, training is performed on a more powerful GPU (NVIDIA A100 with 80 GB VRAM) using a batch size of 32. Training 2D-kFNO and 2D-FNO models requires approximately 43 hours and 72 hours, respectively, with 2D-kFNO benefiting from faster training.

The 1D FNO-based networks are configured with $d_{\z} = 30$ and $\kk_\text{max} = 128$. For 2D FNO-based networks, the configuration is set to $d_{\z} = 20$, $\kk_{\text{1,max}} = \kk_{\text{2,max}} = 64$ and $\alpha$=1. Additionally, to reduce the size of the 2D models, the three Fourier layers within $H$ share the same set of trainable parameters.

All 1D CNN-based networks adopt the following configuration: $L = 7$, $c_1 = 16$, $c_2 = 32$, $c_3 = 64$, $c_4 = 128$, $c_5 = 128$, $c_6 = 64$, and $c_7 = 32$. The kernel size for all convolutional layers is set to 3, and the stride is fixed at 1. The activation function used is ReLU, and layer normalization is applied after each convolutional layer.

\begin{table}
	\caption{ 
		Short-term relative $L_2$ training/validation errors for all models in learning 1D instabilities. Values in parentheses correspond to FNO-based models without skip connections in the Fourier layers. Underlined values indicate models prone to divergence in long-term recurrent predictions.
		\label{Table1} }  
	\centerline{
		\begin{tabular}{ |c | c |c |c | c| }  
			\hline
			     & \multicolumn{2}{c|}{ DL instability: Eq. \eqref{eq:MS} }  & \multicolumn{2}{c|}{ DT instability: Eq. \eqref{eq:KS} }    \\
		          &  $\beta=40$   & $\beta=10$ & $\beta=40$ &$\beta=10$   		  \\
		       & 	Train L2/Valid L2   & 	Train L2/Valid. L2  & 	Train L2/Valid. L2  & 	Train L2/Valid. L2  \\
			\hline
			kFNO  &  \underline{0.0045/0.0063} (0.024/0.025)    & 0.0010/0.0011  &  0.00074/0.00075 & 0.00078/0.00082  \\ 
			kFNO$^\dagger$ &  \underline{0.0036/0.0051} (0.0063/0.030 )      & 0.00057/0.00058  &  0.00059/0.00061  & 0.00075/0.00078 \\ 
			kFNO* &    \underline{0.014/0.022}  (\underline{0.029/0.030})  & 0.0025/0.0025  & 0.0032/0.0034   & 0.0021/0.0022  \\ 
		   FNO  & 0.028/0.028 ( 0.031/0.031 )  & 0.0033/ 0.0034   &  0.0024/0.0025  &  0.0022/0.0023 \\ 
			\hline
			kCNN &  0.026/0.030    & 0.0087/0.0090  & 0.0059/0.0061 & 0.011/0.012\\ 
		    CNN   &  0.039/0.040   & 0.029/0.036 & 0.021/0.021 & 0.018/0.019\\ 
		  \hline
		\end{tabular}
	}
\end{table}

\bibliography{KoopmanLearning}


\newpage
\pagebreak
\pagestyle{plain}

\begin{center}
{\huge Supplementary Material} 
\end{center}

\setcounter{figure}{0}
\setcounter{page}{1}
\makeatletter
\renewcommand{\thefigure}{S\arabic{figure}}

\pagebreak

\begin{figure*}
	\centerline{
		\includegraphics[width=1\linewidth]{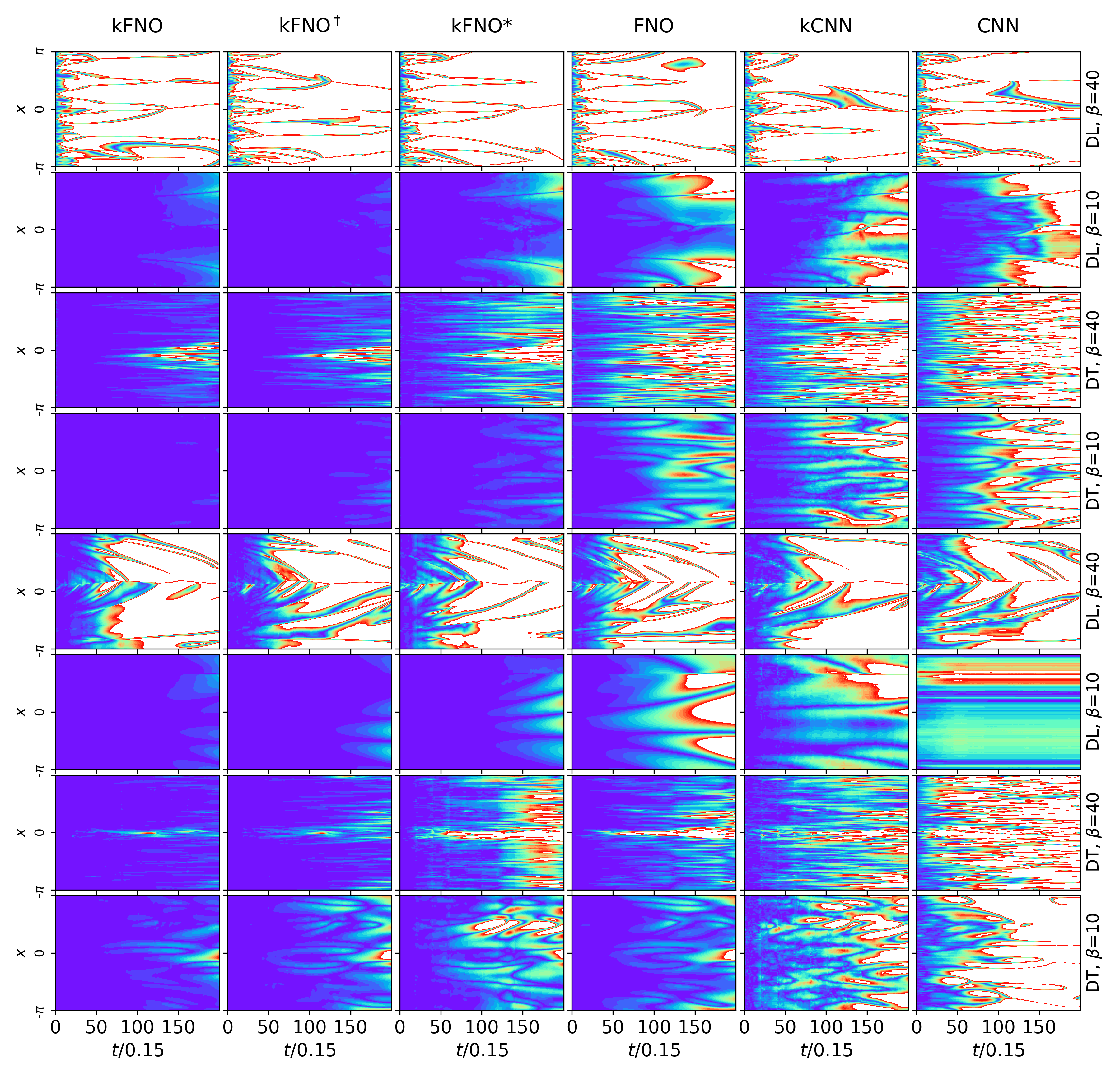}
   }
	\caption{
		\label{fig:err_1D}
		Time evolution of relative $L_2$-error between reference 1D solutions and predictions by kFNO,  kFNO$^\dagger$, kFNO*, FNO, kCNN, and CNN (left to right columns) for four cases (DL/DT at $\beta = 40$ and $\beta = 10$ in rows). The first four rows use random initial conditions, and the last four rows start from long-term solutions evolved after 2000 time steps of size 0.15. Error values from 0 to 0.1 are shown in rainbow colors (blue to red), with values above 0.1 truncated to white.
			}
\end{figure*}

\begin{figure*}
	\centerline{
		\includegraphics[width=1\linewidth]{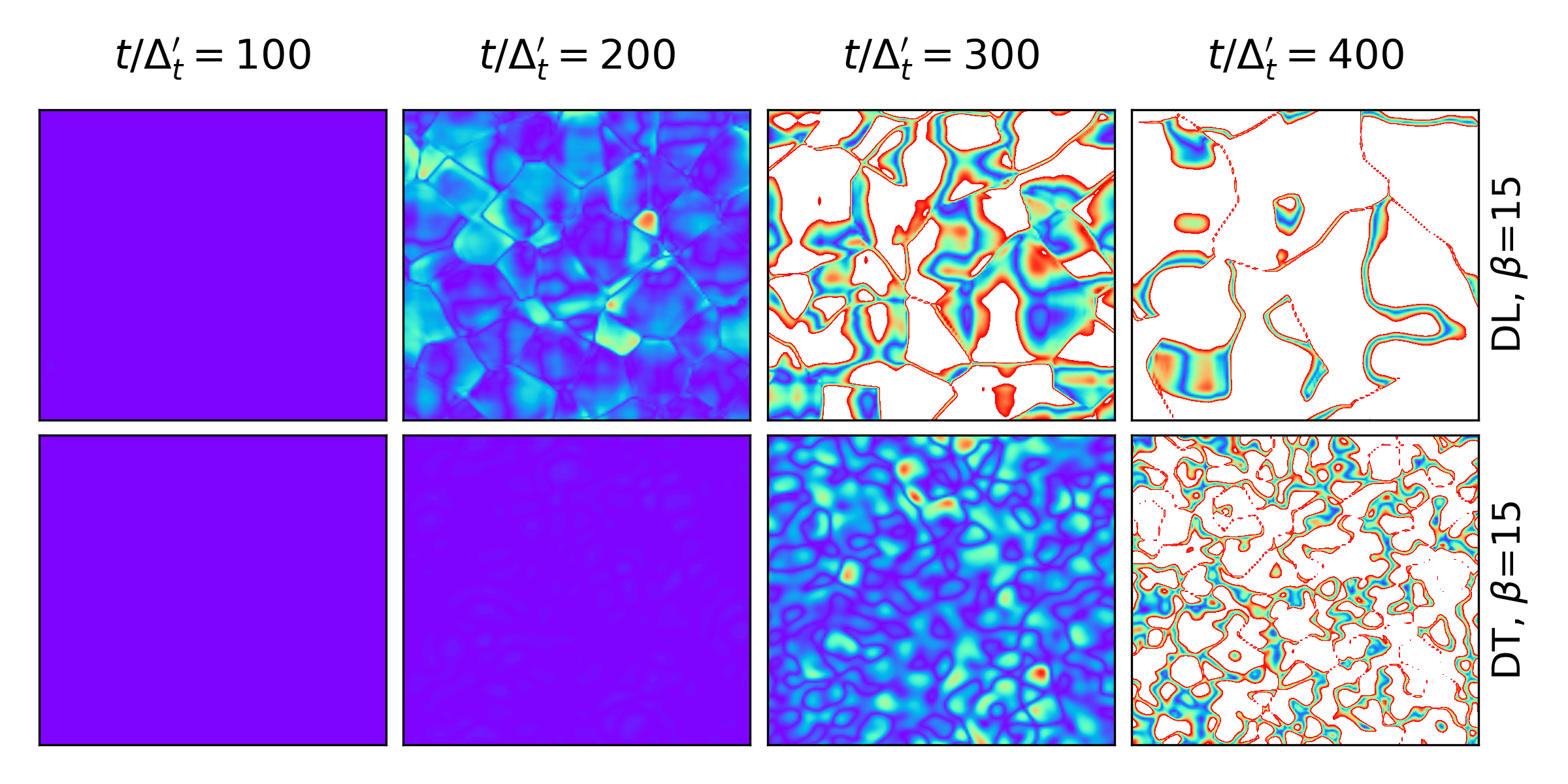	}
   }
	\caption{
		\label{fig:err_2D}
		Temporal development of relative $L_2$ error between reference 2D solutions and kFNO predictions for four cases (DL/DT with $\beta = 15$ in rows), starting from random initial conditions. 
		Errors at $t/\Delta_t' = 100, 200, 300, 400$ are shown from left to right. Error values from 0 to 0.1 are shown in rainbow colors (blue to red), with higher values in white.
		}
\end{figure*}

\end{document}